\documentclass[reqno]{jgokova}
\usepackage{epsfig,graphicx}

\begin{document}

\newcommand{\ben}{\begin{enumerate}}
\newcommand{\een}{\end{enumerate}}
\newcommand{\be}{\begin{equation}}
\newcommand{\ee}{\end{equation}}
\newcommand{\bea}{\begin{eqnarray}}
\newcommand{\eea}{\end{eqnarray}}
\newcommand{\bc}{\begin{center}}
\newcommand{\ec}{\end{center}}

\newtheorem{thm}{Theorem}[section]
\newtheorem{cor}[thm]{Corollary}
\newtheorem{lem}[thm]{Lemma}
\newtheorem{prop}[thm]{Proposition}
\newtheorem{conj}[thm]{Conjecture}

\theoremstyle{definition}
\newtheorem{defn}[thm]{Definition}

\theoremstyle{remark}
\newtheorem{rem}[thm]{\rm\bfseries{Remark}}
\newtheorem*{notation}{Notation}

\newtheorem{ques}[thm]{\rm\bfseries{Question}}
\newtheorem{cons}[thm]{\rm\bfseries{Construction}}
\newtheorem{exm}[thm]{\rm\bfseries{Example}}

\def\C{{\mathbb C}}
\def\P{{\mathbb P}}
\setcounter{page}{1}
\volume{3}

\title[]{The Catanese-Ciliberto-Mendes Lopes surface}
\author[]{Selman Akbulut}

\thanks{The author is partially supported by NSF grant DMS 0905917}


\address{Department  of Mathematics, Michigan State University,  MI, 48824}
\email{akbulut@math.msu.edu}

\begin{abstract}
We draw a handlebody picture of the complex surface defined by  Catanese-Ciliberto-Mendes Lopes. This is a surface obtained by taking the quotient  of the product of surfaces  $\Sigma_{2} \times \Sigma_{3}$  of genus $2$ and $3$, under the product of involutions $\tau_{2} \times \tau_{3}$, where $\tau_{2}$ is the elliptic involution of $\Sigma_{2}$, and $\tau_{3}$ is a free involution on $\Sigma_{3}$.
\end{abstract}

\maketitle

\section{Introduction}

\vspace{.1in}

{\it Catanese-Ciliberto-Mendes Lopes} surface (CaCiMe or CCM surface in short) $M$ is a  complex surface constructed in \cite{ccm} (and discussed in \cite{hp} and \cite{p}), which is  topologically a genus $2$ surface bundle over a surface of genus $2$. Recently in \cite{ap}  this surface is used in  interesting smooth manifold constructions. While inspecting \cite{ap} we felt that first  this interesting complex surface $M$ deserves a careful topological study of its own. Generally, drawing  handlebody pictures of circle bundles over $3$-manifolds, or $3$-manifold bundles over circles is relatively easy compare to surface bundles of surfaces (e.g.\cite{ak} and \cite{a1}). Here we take this opportunity to introduce a new technique to draw a surface bundle over a surface, which avoids ``turning handles upside down'' process.   This manifold $M$ is a good test case to understand many of the general difficulties one encounters in drawing surface bundles over surfaces, as well as taking their fiber sum. We draw the handlebody picture of $M$ in such a way that all the tori used in the constructions of \cite{ap} are clearly visible. This combined with  the log transform picture (e.g. \cite{ay}) will allow one to see the ``Lutinger surgery'' constructions of \cite{ap}  in a concrete geometric way.

\vspace{.05in}

\section{Construction}

\vspace{.05in}

Let $\Sigma_{g}$ denote the surface of genus $g$. Let   $\tau_{2} :\Sigma_{2}\to \Sigma_{2}$ be the hyperelliptic involution and $\tau_{3}: \Sigma_{3}\to \Sigma_{3}$ be the free involution  induced by $180^{o}$ rotation, as indicated in Figure 1.  The CaCiMe surface $M$  is the complex surface obtained by taking the quotient of $\Sigma_{2} \times \Sigma_{3}$ by the product involution: 

$$M=(\Sigma_{2} \times \Sigma_{3} ) / \tau_{2}\times \tau_{3}$$

 \begin{figure}[ht]  \begin{center}  
\includegraphics[width=.4\textwidth]{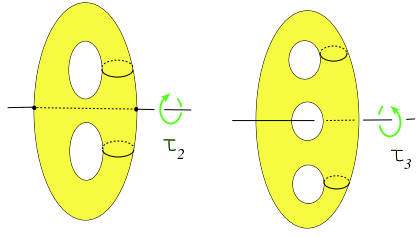}   
\caption{Involutions $\tau_{2}$ and $\tau_{3}$} 
\end{center}
\end{figure} 

\vspace{.1in}

By projecting to the second factor we can describe  $M$ as a  $\Sigma_{2}$-bundle over $\Sigma_{2}=\Sigma_{3}/\tau_{3}$. Let $A$ denote the twice punctured  $2$-torus $A=T^{2}-D^{2}_{-}\sqcup D^2_{+}$. Then clearly $M$ is obtained by identifying the two boundary components of  $\Sigma_{2} \times A$ by the involution induced by $\tau_{2}$ (notice $A$ is the interior of the fundamental domain of the action $\tau_{3}$). By  deforming $A$ as in Figure 2, we see that $M= E\;\natural \;E'$  is obtained by fiber summing two $\Sigma_{2}$ bundles over $T^{2}$,  where $E$ is the trivial bundle $\Sigma_{2}\times T^{2}\to T^{2}$, and 

$$E'= \Sigma_{2}\times S^1 \times [0,1] / (x,y,0)\sim (\tau_{2}(x),y,1) $$

 \begin{figure}[ht]  \begin{center}  
\includegraphics[width=.75\textwidth]{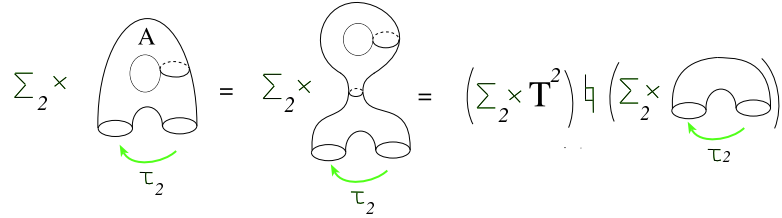}   
\caption{Decomposing $M= E \; \natural \;E'$} 
\end{center}
\end{figure} 

We will build a handlebody of $M$ by a step by step process drawing the following handlebodies in the given order,  also we will concretely identify the indicated diffeomorphisms:

\begin{itemize}

\vspace{.1in}

\item[(a)]  $E_{0} =E - \Sigma_{2} \times D^{2} =\Sigma_{2} \times (T^2 - D^2)$ 
 
\item[(b)]  $E'_{0} = E'-  \Sigma_{2} \times D^2$ 

\item[(c)] $ f_{1}: \partial E_{0} \stackrel{\approx}{\longrightarrow} \Sigma_{2} \times S^1 $ 

\item[(d)]  $ f_{2}: \partial E'_{0} \stackrel{\approx}{\longrightarrow} \Sigma_{2} \times S^1 $ 

\item[(e)] $M= -E_{0} \smile_{ f_{2}^{-1}\circ f_{1}}  E'_{0}$

\end{itemize}

\vspace{.1in}

One way to to perform the gluing operation (e) is to turn the handlebody $E_0$ upside down and attach its dual handlebody to top of $E_{0}'$, getting $M=-E_{0} \smile E_{0}'$ (e.g. the  technique used in \cite{a1}). In this paper we choose another way which amounts to identifying the boundaries of the two handlebodies $-E_{0}$ and $E_{0}'$ by a cylinder $\partial E_{0} \times I $ $$M= -E_{0} \smile_{f_{1}} (\Sigma_{2}\times S^1)\times I\smile_{f_{2}^{-1}} E_{0}'$$
Though this seems a trivial distinction, it makes a big difference in constructing the handlebodies. One advantage of this technique is that we see the imbedded tori used in the construction of \cite{ap} clearly.

\vspace{.05in}

\section{Constructing $E_{0}$}

\vspace{.05in}

Figure 3 (a disk with two pairs of $1$-handles and a $2$-handle, where only the attaching arcs of the $1$-handles are drawn) describes a handlebody for $\Sigma_{2}$, and  Figure 4 is  $\Sigma_{2}\times [0,1]$. Hence Figure 5 is a handlebody of $\Sigma_{2}\times S^1$ (compare \cite{ak}). Figure 6 is the same as Figure 5, except it is drawn as a Heegard diagram.  So Figure 7 describes a handlebody picture of of $\Sigma_{2}\times S^1 \times S^1=\Sigma_{2}\times T^2$. A close inspection shows that removing the $2$-handle, denoted by $c$, from Figure 7 gives the handlebody of $E_{0}$ ($c$ is the disk boundary in $T^{2}- D^2$, as the attaching circle of the $2$-handle corresponding to $D^2$; more precisely it is the upside down $2$-handle of the missing $\Sigma_{2}\times D^2$ which is removed from $E_{0}$).  

\vspace{.05in}

Next in  Figures 8 through 11 we gradually convert the  ``pair of balls'' notation of the $1$-handles of Figure 7 to the ``circle-with-dot" notation of  \cite{a2} (i.e. carving). Figure 11 is the same as Figure 7, except that all of its $1$-handles are drawn in circle-with-dot notation. For the benefit of the reader we did this transition in several steps:  First in Figure 8 we converted a pair of $1$-handles of Figure 7 to the circle-with-dot notation, then in Figure 11 converted the remaining $1$-handles. Figure 9 shows how to perform local isotopies near the attaching balls of $1$-handles to go to intermediate picture Figure 10 were the attaching balls are drawn as flat arcs. We then converted the flat arcs to the circle-with dots. In Figure 11 all  the 2-handles are attached with $0$-framing.

\vspace{.05in}

\section{Diffeomorphism $f_{1}: \partial E_{0} \to \Sigma_{2}\times S^1$}

\vspace{.05in}

Next we construct a difeomorphism $f_{1}: \partial E_{0} \stackrel{\approx}{\longrightarrow} \Sigma_{2} \times S^1$. 
First by an isotopy we go from Figure 11 to Figure 12, then by replacing the
circles with dots with 0-framed circles, and by performing the handle
slides to Figure 12 as indicated by the arrows, we obtain the first picture of Figure 13, and by further handle slides and cancellations we obtain the second picture of Figure 13, which is $\Sigma_{2}\times D^2$. In the Figure 13 we also indicate where this diffeomorphism throws the linking loops $a_1, b_1, a_2, b_2, c$. Finally in Figure 14 we describe the diffeomprphism we constructed $f_{1}: \partial E_{0} \to \Sigma_{2}\times S^1$ in a much more concrete way by indicating the images of the arcs shown in the figure. Though going from Figure 11 to Figure 14 is locally a routine process, finding the correct handle sliding moves and locating and keeping the track of those arcs is the most time consuming part of this work.

\vspace{.05in}

\section{Constructing $E_{0}'$ and  $f_{2}:\partial E_{0}'\to \Sigma_{2}\times S^1$}

\vspace{.05in}

Figures 15 and 16 shows that the diffoemorphism $\tau_{2}: \Sigma_{2}\to \Sigma_{2}$ is induced from $180^{0}$ rotation of the disk with four $1$-handles. Having noted this, we proceed exactly as Figure 7 through Figure 11, except that we replace Figure 7 by Figure 17 (due to twisting by $\tau_{2}:\Sigma_{2}\to \Sigma_{2}$). So Figure 20 is the handlebody of $E_{0}'$ (without the curve denoted by $c'$), and Figure 22 describes a diffeomorphism $f_{2}:\partial E_{0}'\to \Sigma_{2}\times S^1$

\vspace{.05in}

\section{Constructing $M=-E_{0}\smile E_{0}'$}

\vspace{.05in}

To construct  $M=-E_{0}\smile E_{0}'$ we draw the  handlebodies $-E_{0}$ and $E_{0}'$ side by side, and glue their boundaries to the two boundary components of the cylinder $\Sigma_{2}\times S^1 \times I$. This gluing is done by identfying  $1$-handle circles $\{a_1, b_1, a_2, b_2, c \}$  (Figure 13) of $\partial E_{0} \approx  \Sigma_{2}\times S^1 $ and  of $\partial E_{0}'\approx \Sigma_{2}\times S^1$ by 2-handles. Now Figure 14 and Figure 22 gives us exactly the information needed to draw   the CaCiMe surface $M=-E_{0}\smile_{f_{2}^{-1}\circ f_{1}} E_{0}'$ as shown in Figure 23 (all the circles are $0$-framed $2$-handles).

\vspace{.05in}

\section{Epilogue}

\vspace{.05in}

 Cacime surface has its own place in the classification scheme of complex surfaces, as stated in the following theorem:

\begin{thm}(\cite{hp}, \cite{p}) If $X$ is a smooth minimal complex projective surface of general type with $p_{g}(X) = q(X) = 3$, then one of following hold:
 \begin{itemize}
\item [(a)] $K_{X}^{2}=6 $ and $X= Sym^{2}(\Sigma_{3})$
\item [(b)] $K_{X}^{2}=8 $ and $X$= Cacime surface.
\end{itemize}
\end{thm}

Recall  $b_{1}(X)=2q(X)$,  $K_{X}=c_{1}(X)$ and $3\sigma(X)=c_{1}^{2}(X)-2 \chi(X)$,   Noether formula: $1-q(X) +p_{g}(X)=\frac{1}{12} [ \;c_{1}^{2}(X)+c_{2}(X) \;]$ $ \Rightarrow $ So if (b) holds  $b_{2}(X)=14$ ad $\sigma(X)=0$
So the Cacime surface $X$ is homology equivalent to $ \# 7(S^{2}\times S^2)\# 6(S^1\times S^3)$, and its fundamental group presumably can be  calculated from its fibration structure provided its monodromies are determined,  But now we can easily calculate the fundamental group as well as the other topological invariants from its handlebody picture in Figure 23. 

\vspace{.04in}

\begin{rem}
One of the reason we decided to study the handlebody structure of the Cacime surface is that, it appears to be the starting point of many other interesting manifolds, for example the construction techniques used in  \cite{a3}, \cite{a4} and \cite{a5}  are all driven from the construction of the handlebody of the Cacime surface. In particular the reader is encouraged to look at  \cite{a3}, where a simpler version of the construction of Section 3 is used to build handlebodies for $T^4$ and $T^2 \times (T^{2} - D^2)$.
\end{rem}


\noindent {\bf Acknowledgements:} We would like to thank Anar Akhmedov for introducing us to   the Catanese-Ciliberto-Mendes Lopes  complex surface, and explaining \cite{ap}. We also  thank IMBM (Istanbul mathematical sciences research institute) for providing us an inspiring environment where a large part of this work was done.

\newpage

 \begin{figure}[ht]  \begin{center}  
\includegraphics[width=.35\textwidth]{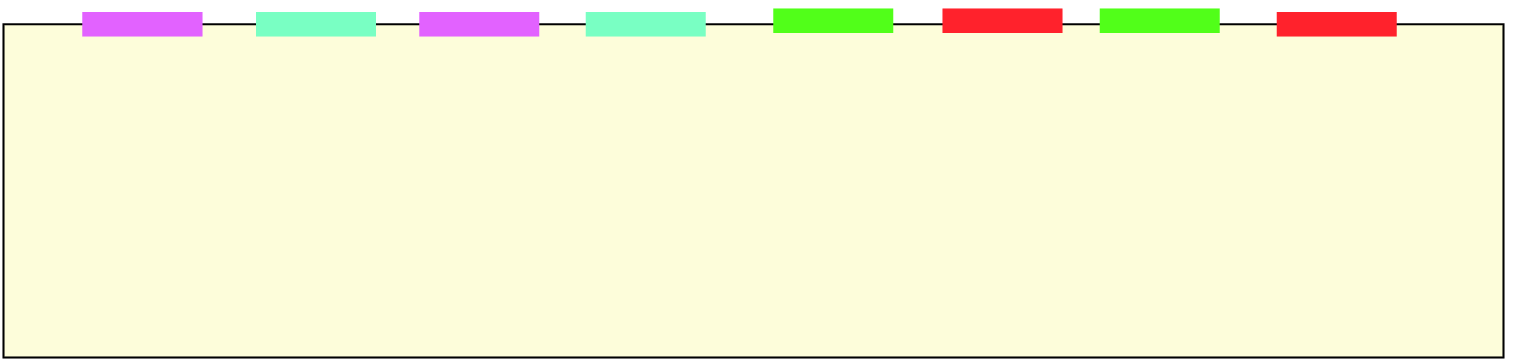}   
\caption{$\Sigma_{2} $} 
\end{center}
\end{figure}

 \begin{figure}[ht]  \begin{center}  
\includegraphics[width=.4\textwidth]{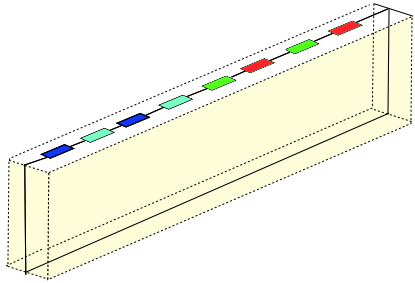}   
\caption{$\Sigma_{2} \times [0,1]$} 
\end{center}
\end{figure}

 \begin{figure}[ht]  \begin{center}  
\includegraphics[width=.8\textwidth]{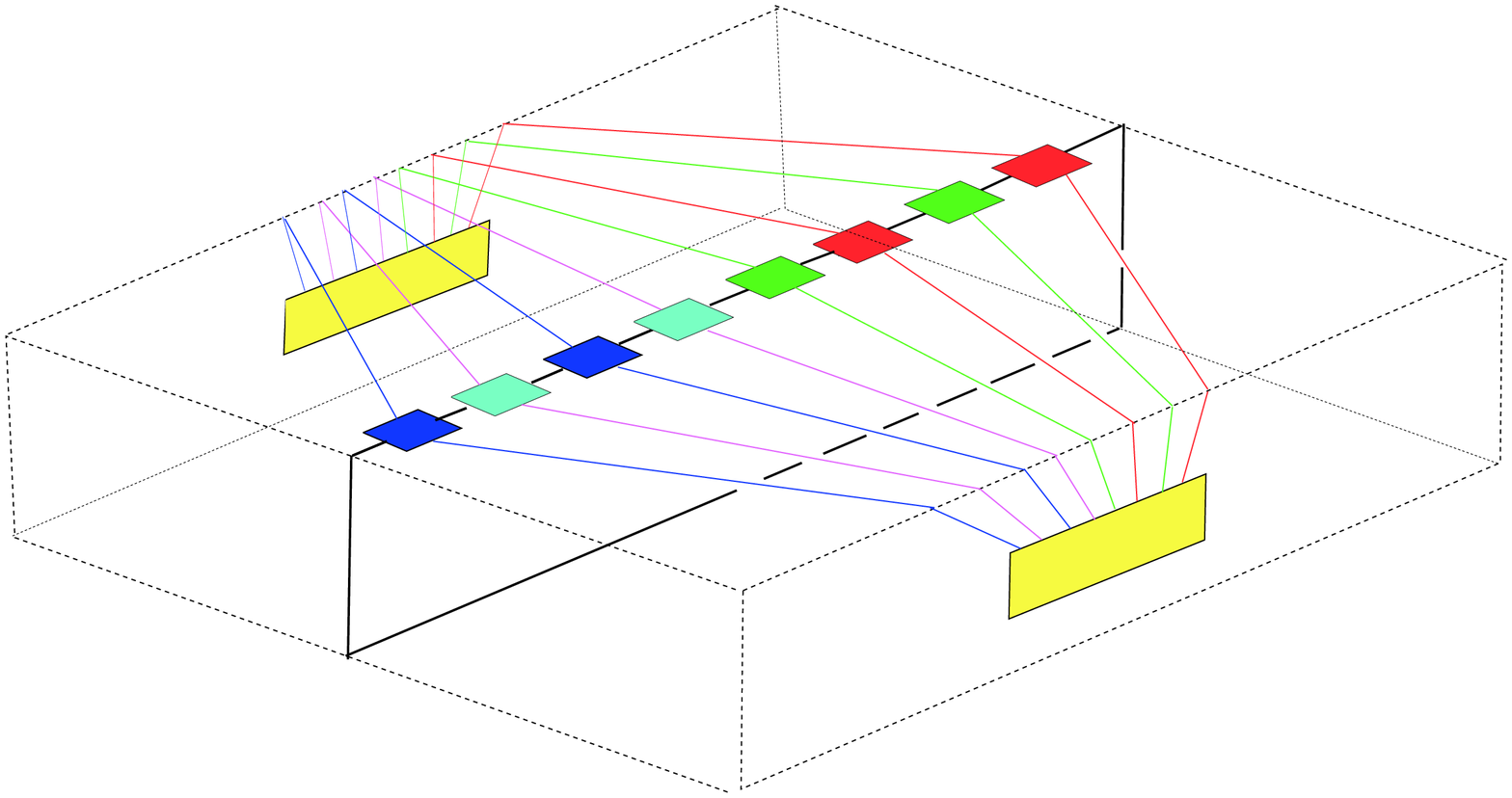}   
\caption{$\Sigma_{2} \times S^1$} 
\end{center}
\end{figure}

 \begin{figure}[ht]  \begin{center}  
\includegraphics[width=.35\textwidth]{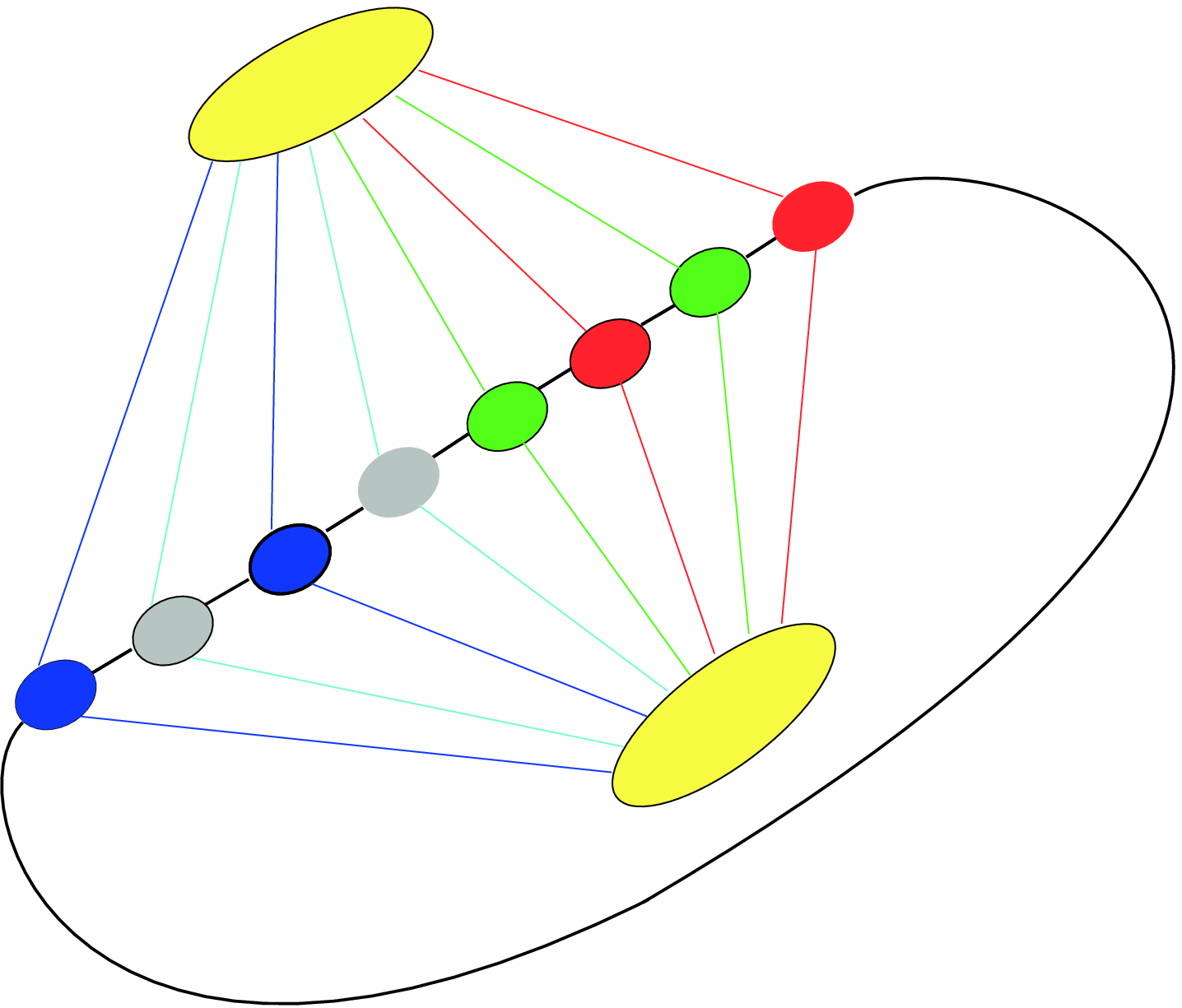}   
\caption{$\Sigma_{2} \times S^1$} 
\end{center}
\end{figure}

\begin{figure}[ht]  \begin{center}  
\includegraphics[width=.6\textwidth]{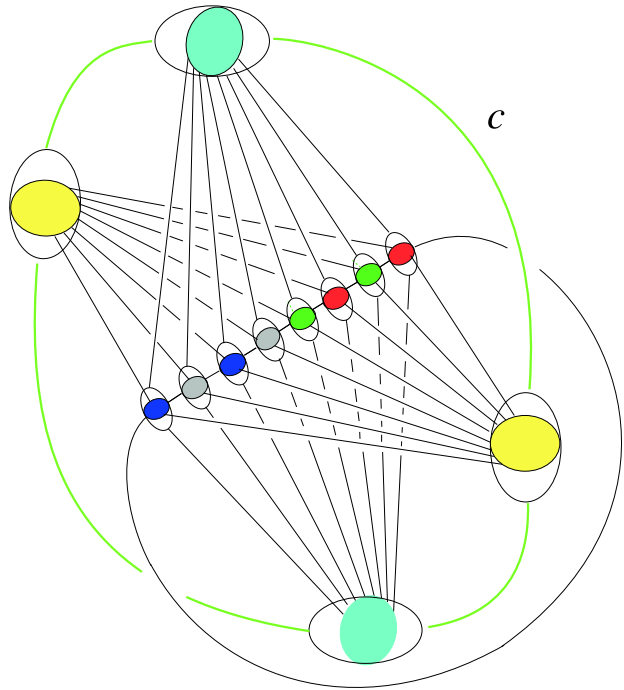}   
\caption{$\Sigma_{2} \times T^2$} 
\end{center}
\end{figure}

 \begin{figure}[ht]  \begin{center}  
\includegraphics[width=.7\textwidth]{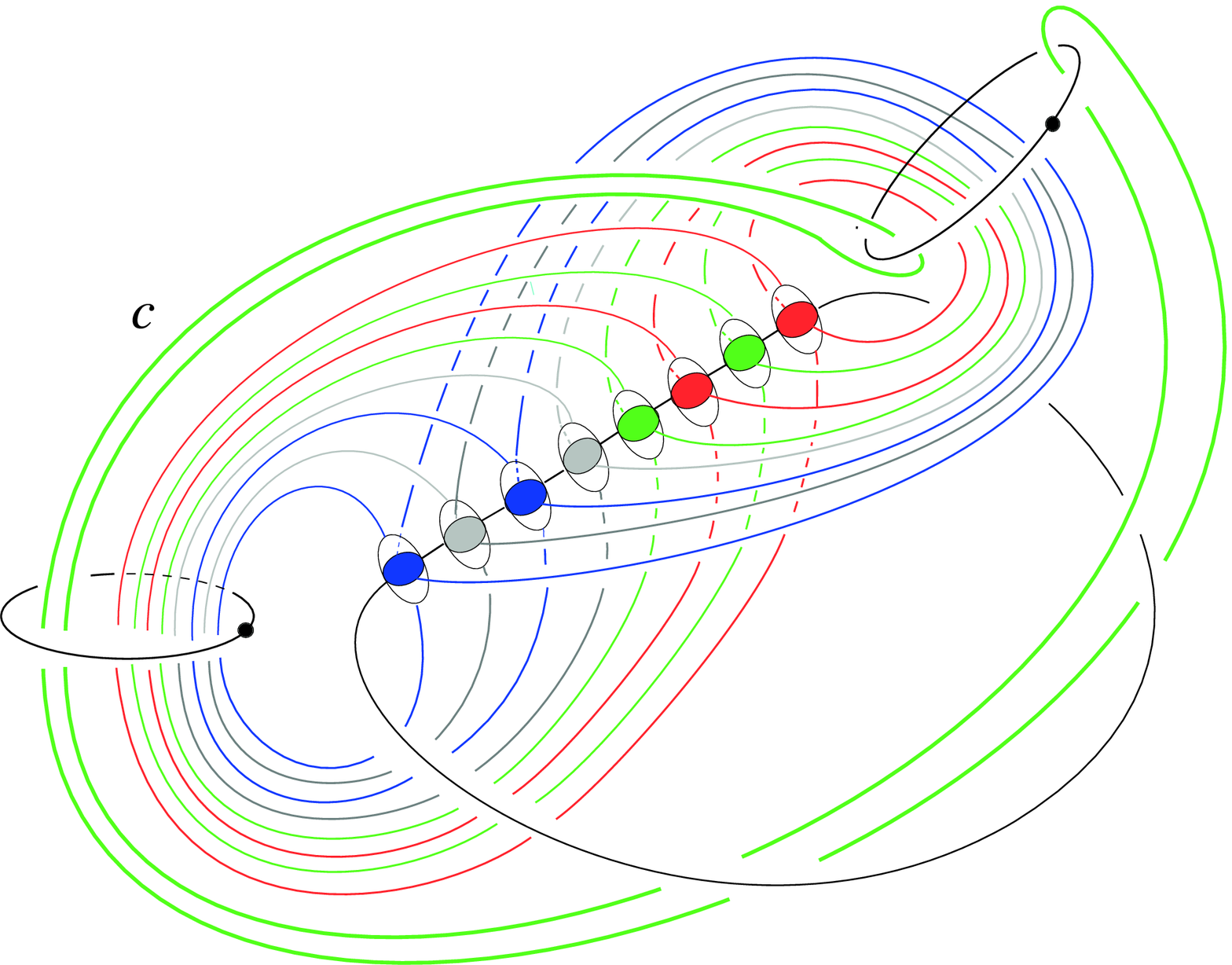}   
\caption{$E_{0}$ } 
\end{center}
\end{figure}

 \begin{figure}[ht]  \begin{center}  
\includegraphics[width=.3\textwidth]{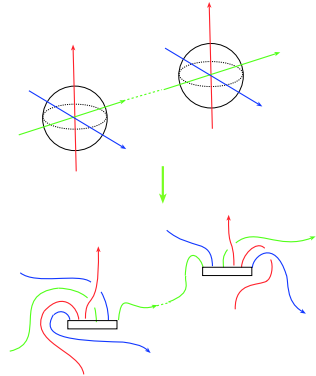}   
\caption{Local isotopies used to flatten $1$-handle attaching balls} 
\end{center}
\end{figure} 

 \begin{figure}[ht]  \begin{center}  
\includegraphics[width=.8\textwidth]{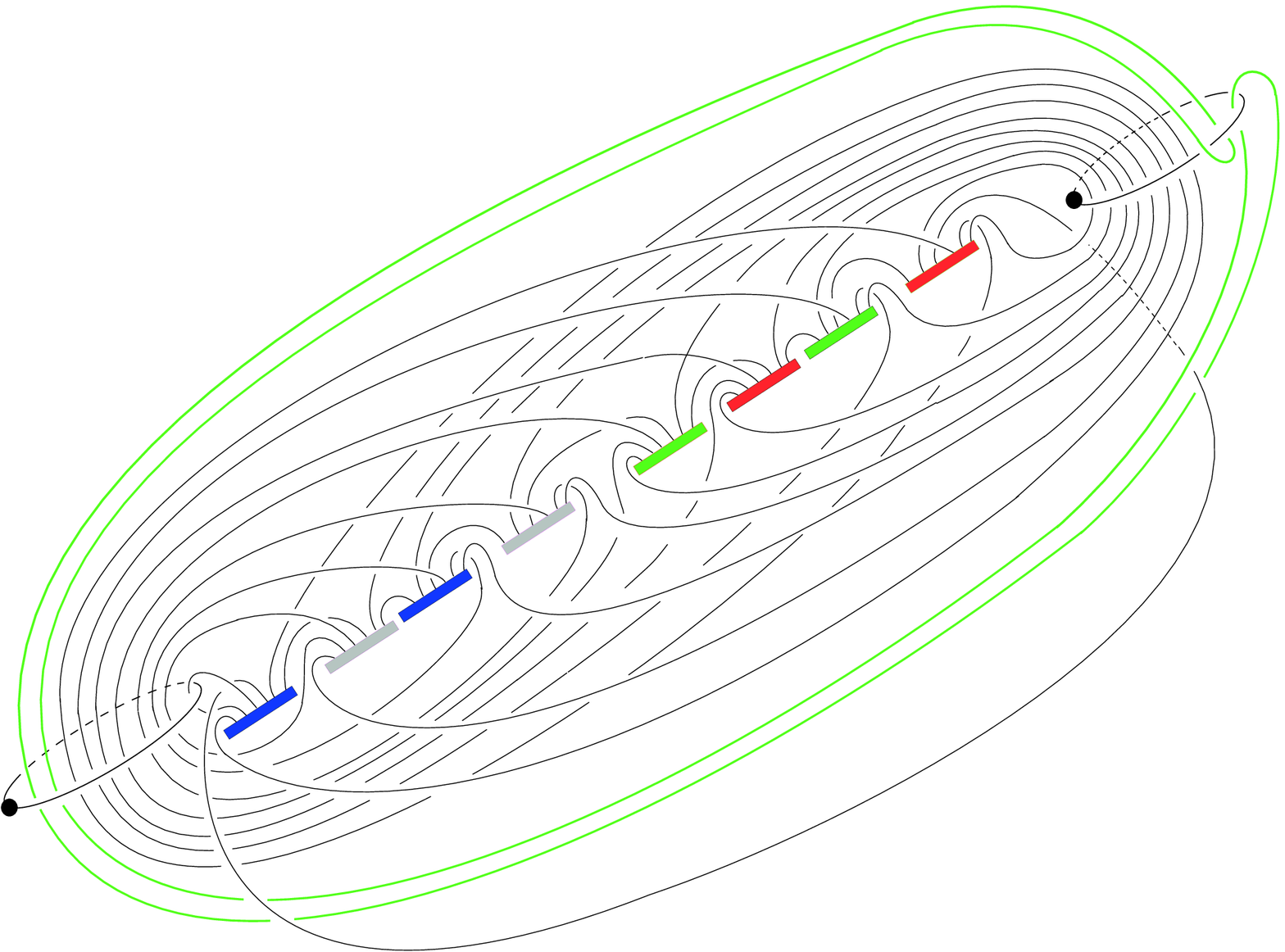}   
\caption{Converting $1$-handle notation} 
\end{center}
\end{figure}

\begin{figure}[ht]  \begin{center}  
\includegraphics[width=.8\textwidth]{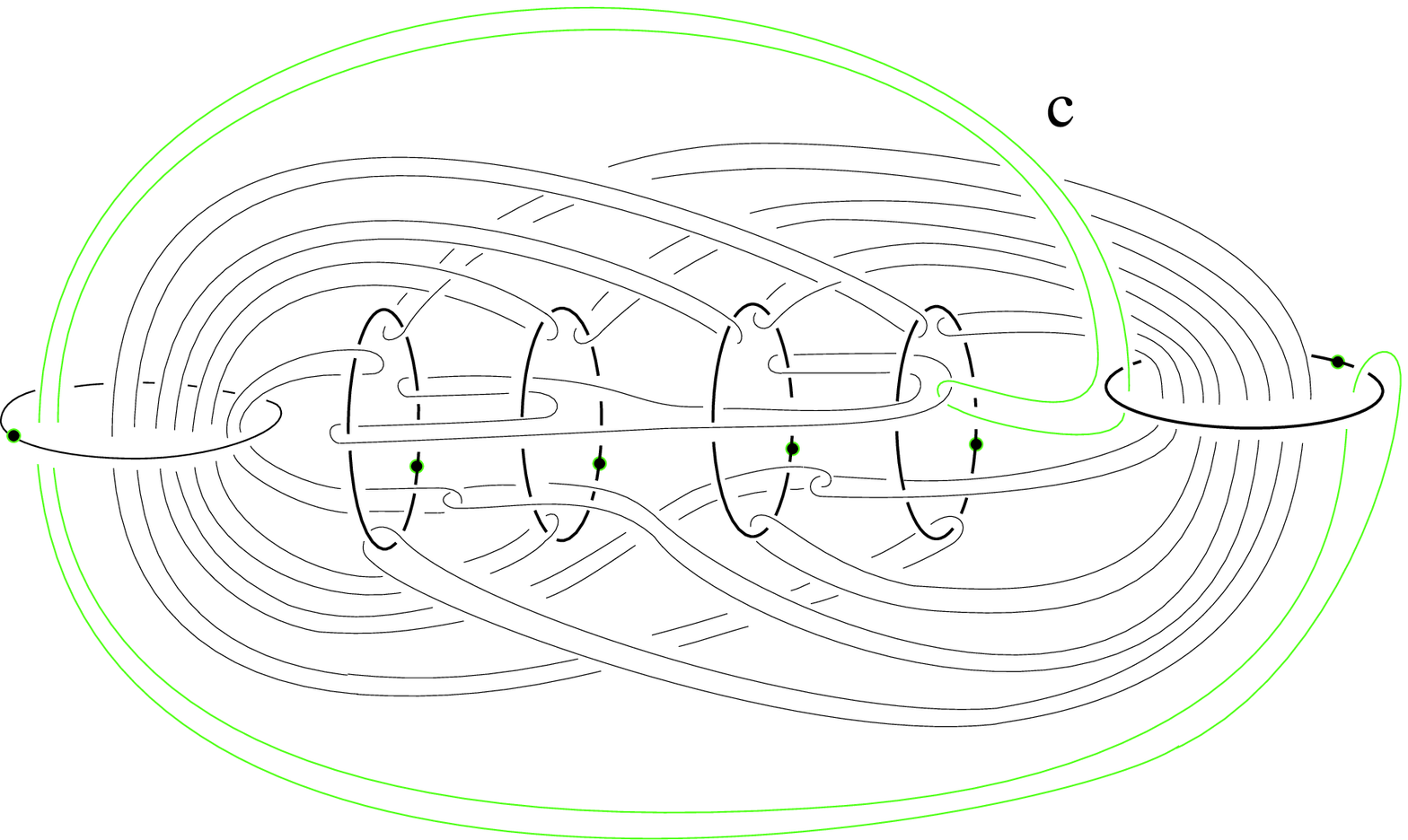}   
\caption{$E_{0}$} 
\end{center}
\end{figure}

\begin{figure}[ht]  \begin{center}  
\includegraphics[width=.7\textwidth]{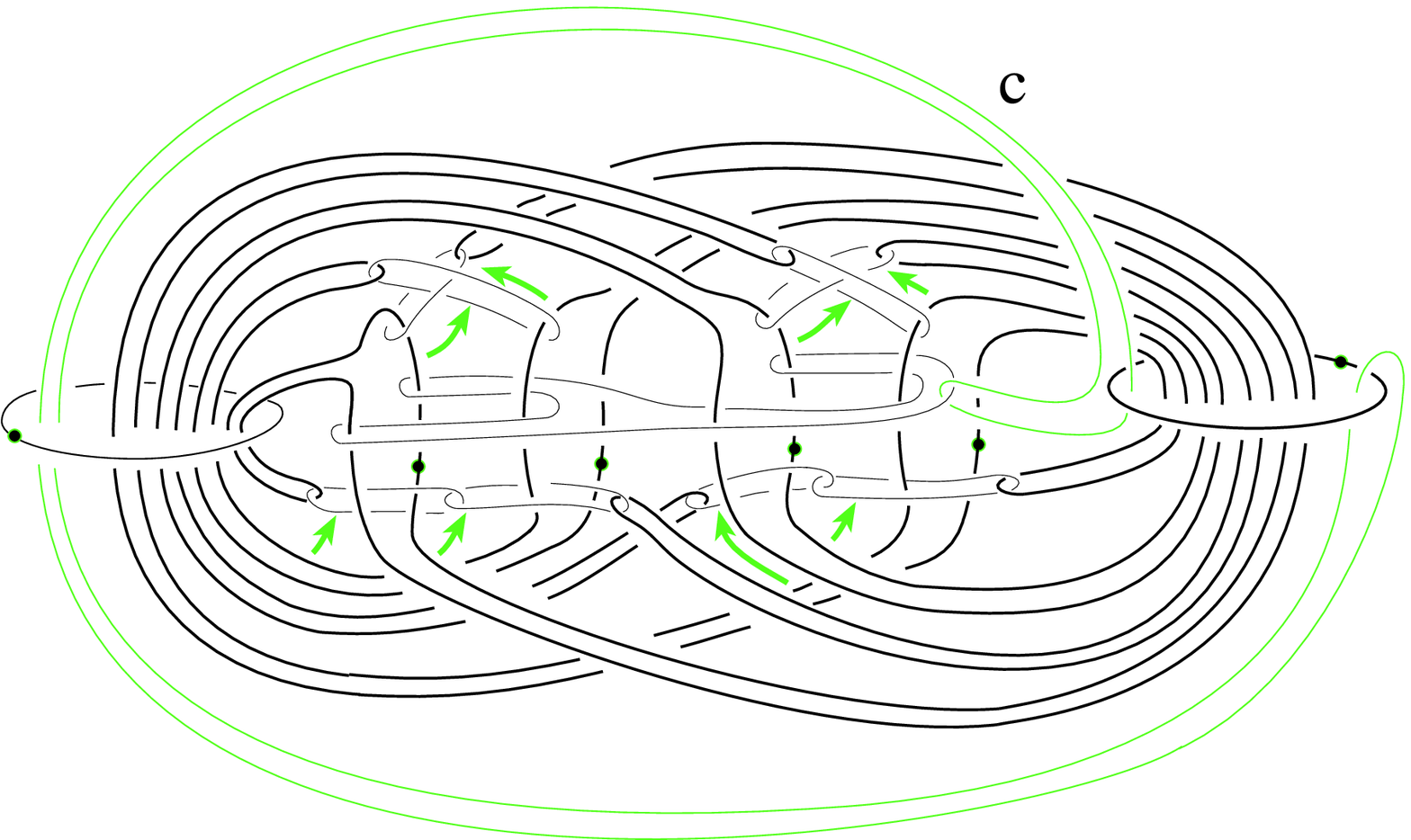}   
\caption{Surgering inside of $E_{0}$} 
\end{center}
\end{figure}

\begin{figure}[ht]  \begin{center}  
\includegraphics[width=.65\textwidth]{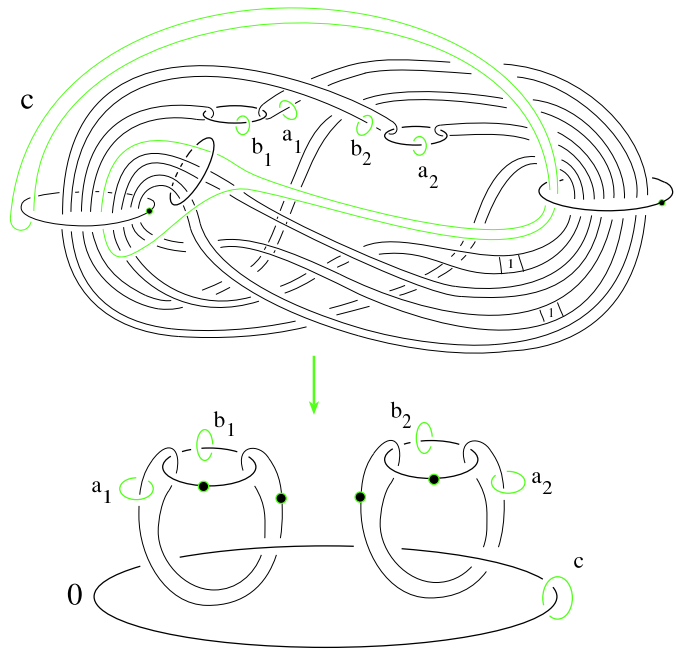}   
\caption{Checking $\partial E_{0}\approx \Sigma_{2}\times S^1$} 
\end{center}
\end{figure}

\newpage

\begin{figure}[ht]  \begin{center}  
\includegraphics[width=.8\textwidth]{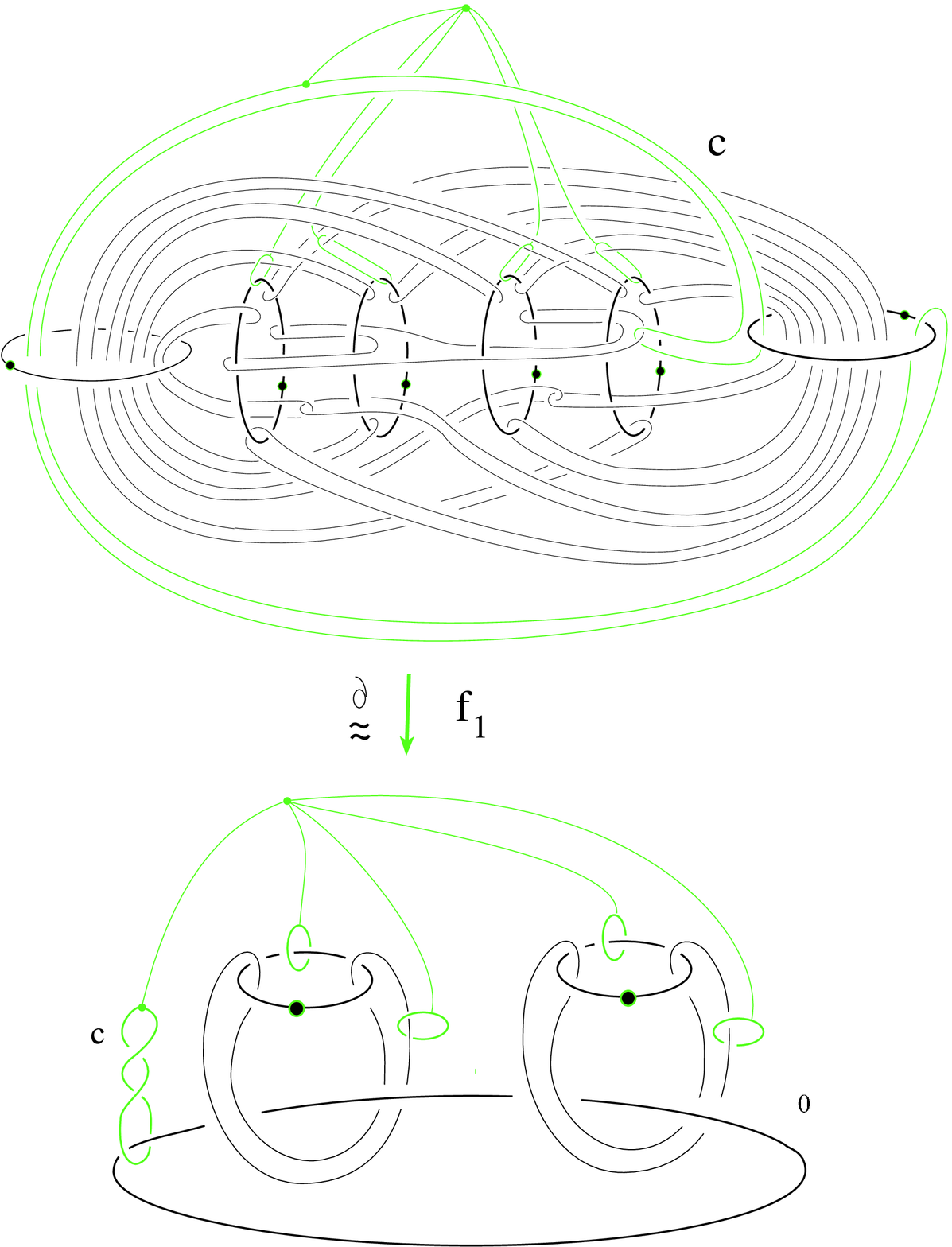}   
\caption{Diffeomorphism $\partial E_{0}\approx \Sigma_{2}\times S^1$ made concrete} 
\end{center}
\end{figure}

 \begin{figure}[ht]  \begin{center}  
\includegraphics[width=.65\textwidth]{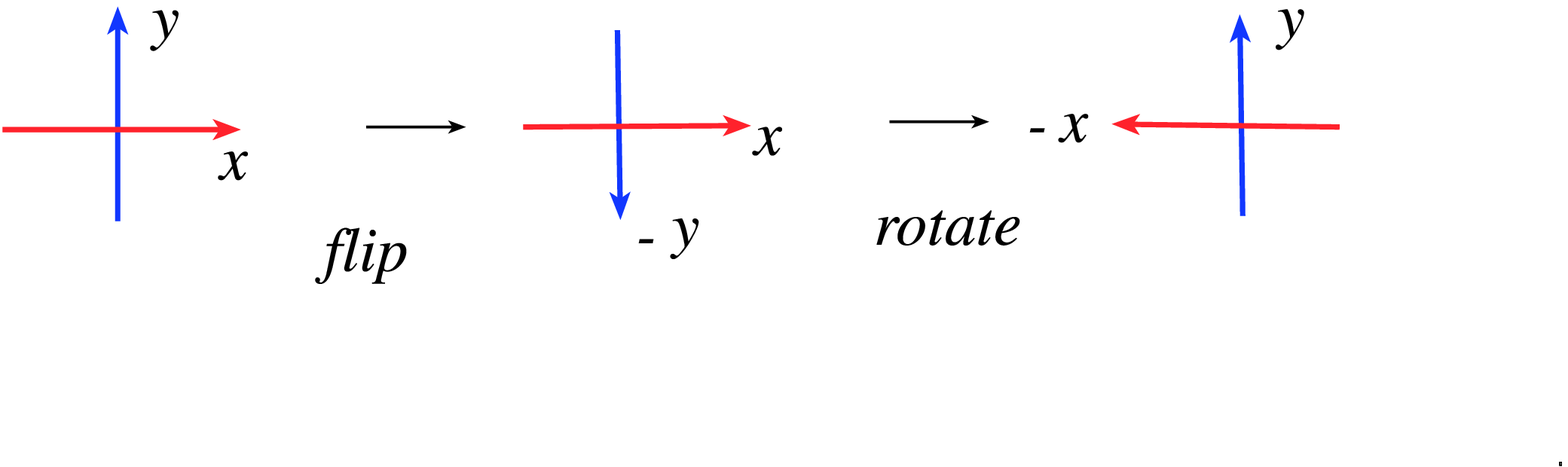}   
\caption{Action of $\tau_{2}$ on the zero handle of $\Sigma_{2}$} 
\end{center}
\end{figure}

 \begin{figure}[ht]  \begin{center}  
\includegraphics[width=.7\textwidth]{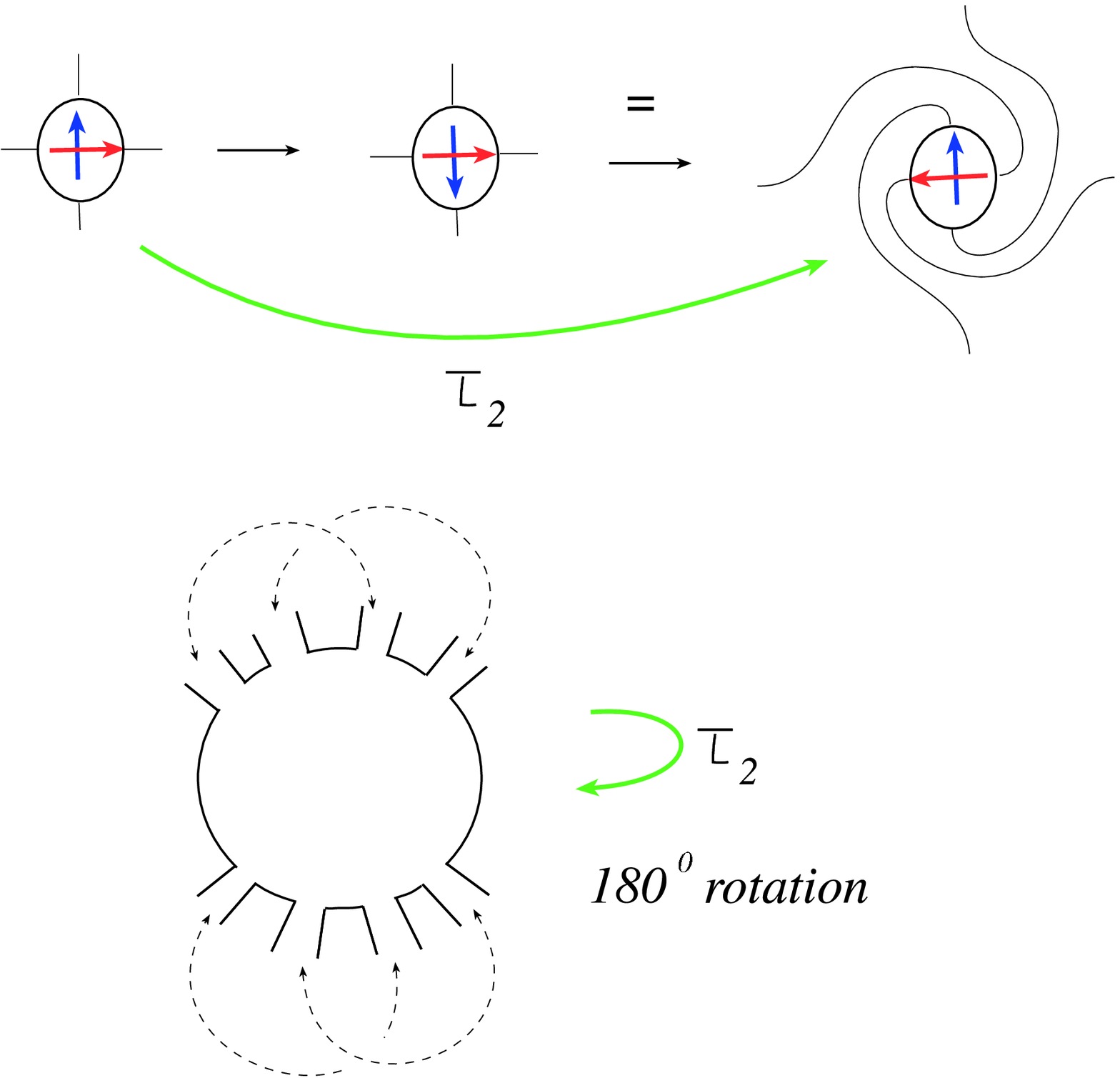}   
\caption{Describing  $\tau_{2}: \Sigma_{2} \to \Sigma_{2}$ } 
\end{center}
\end{figure}

 \begin{figure}[ht]  \begin{center}  
\includegraphics[width=.7\textwidth]{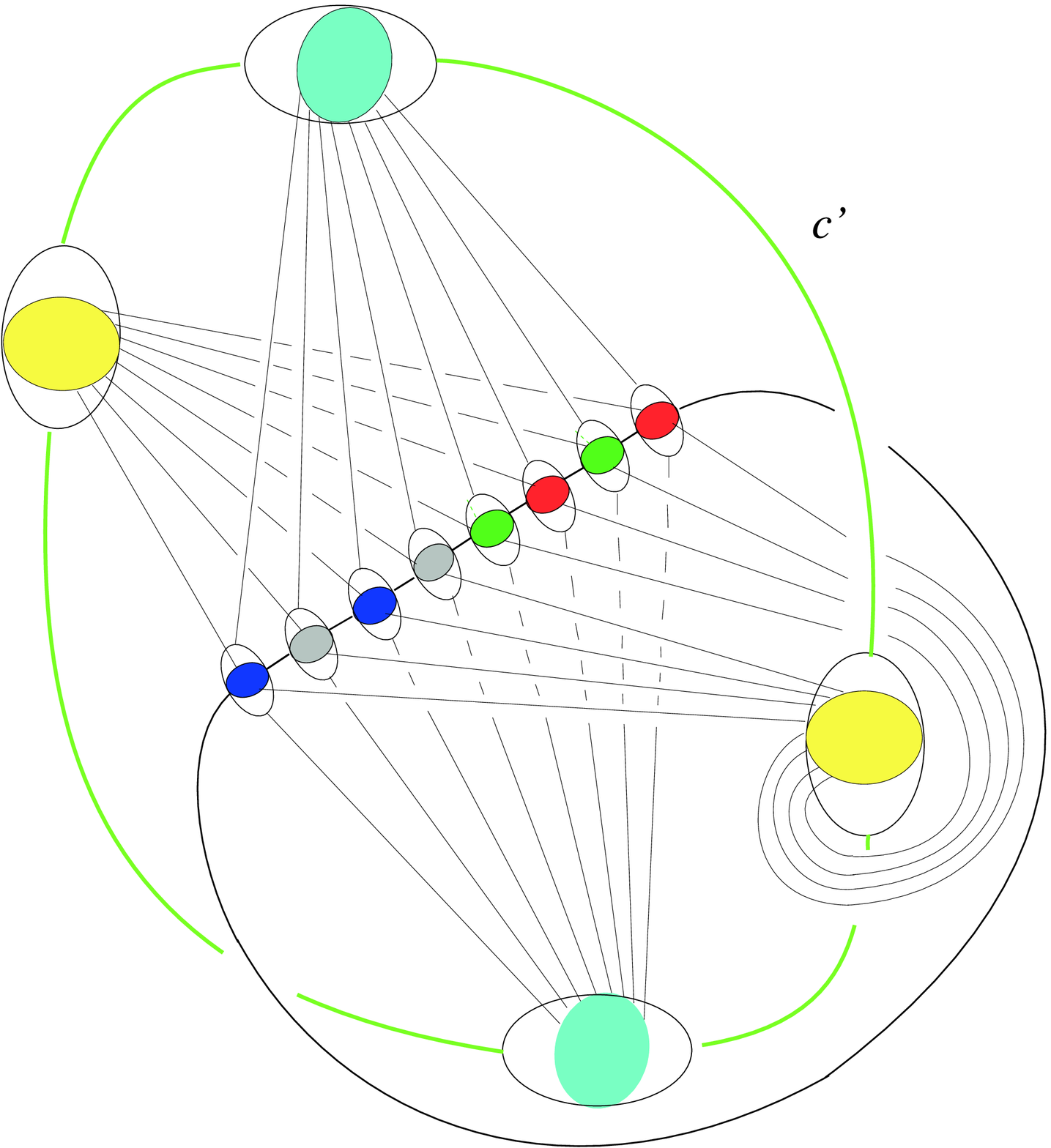}   
\caption{$E_{0}'$} 
\end{center}
\end{figure}

 \begin{figure}[ht]  \begin{center}  
\includegraphics[width=.7\textwidth]{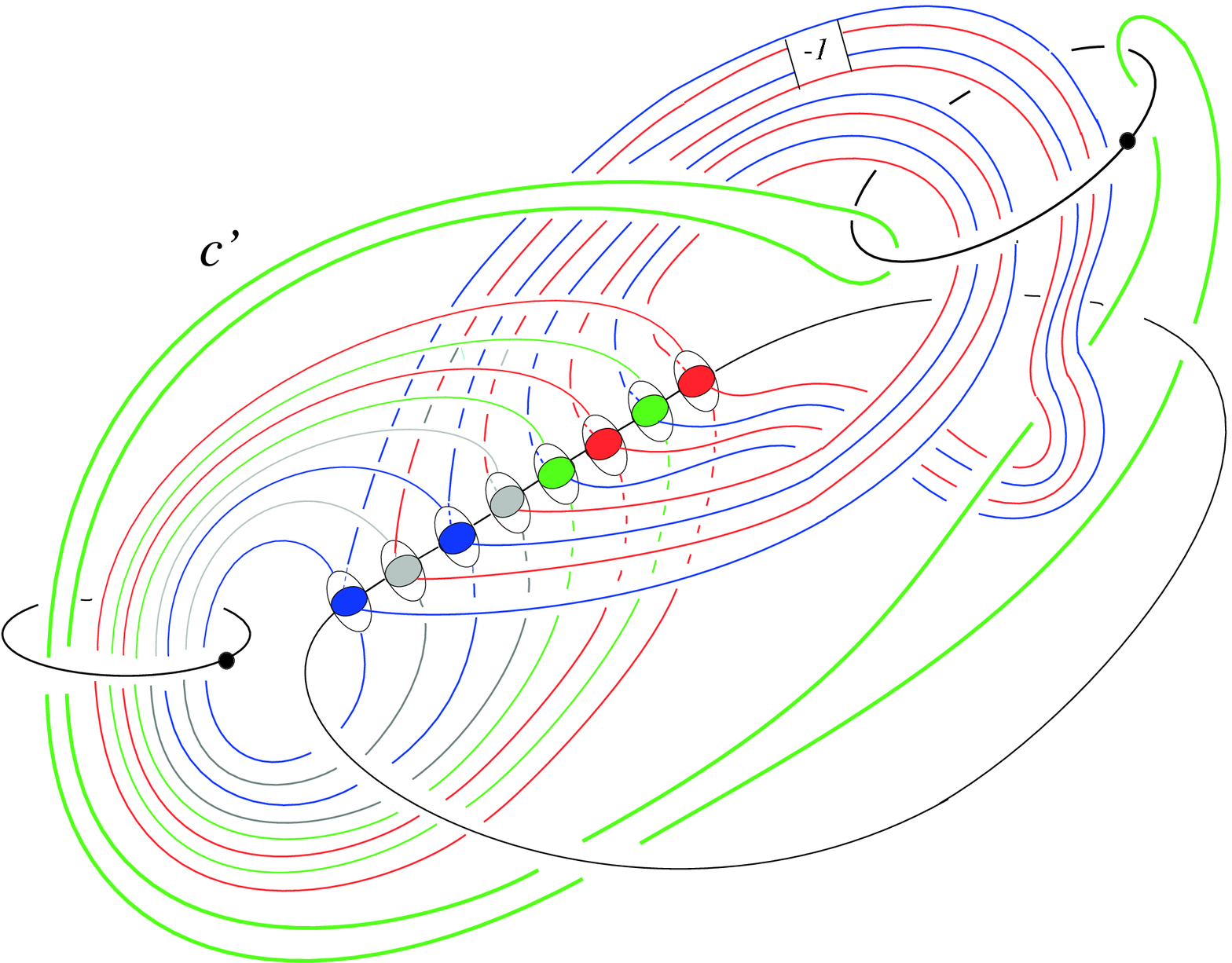}   
\caption{$E_{0}'$} 
\end{center}
\end{figure}

 \begin{figure}[ht]  \begin{center}  
\includegraphics[width=.8\textwidth]{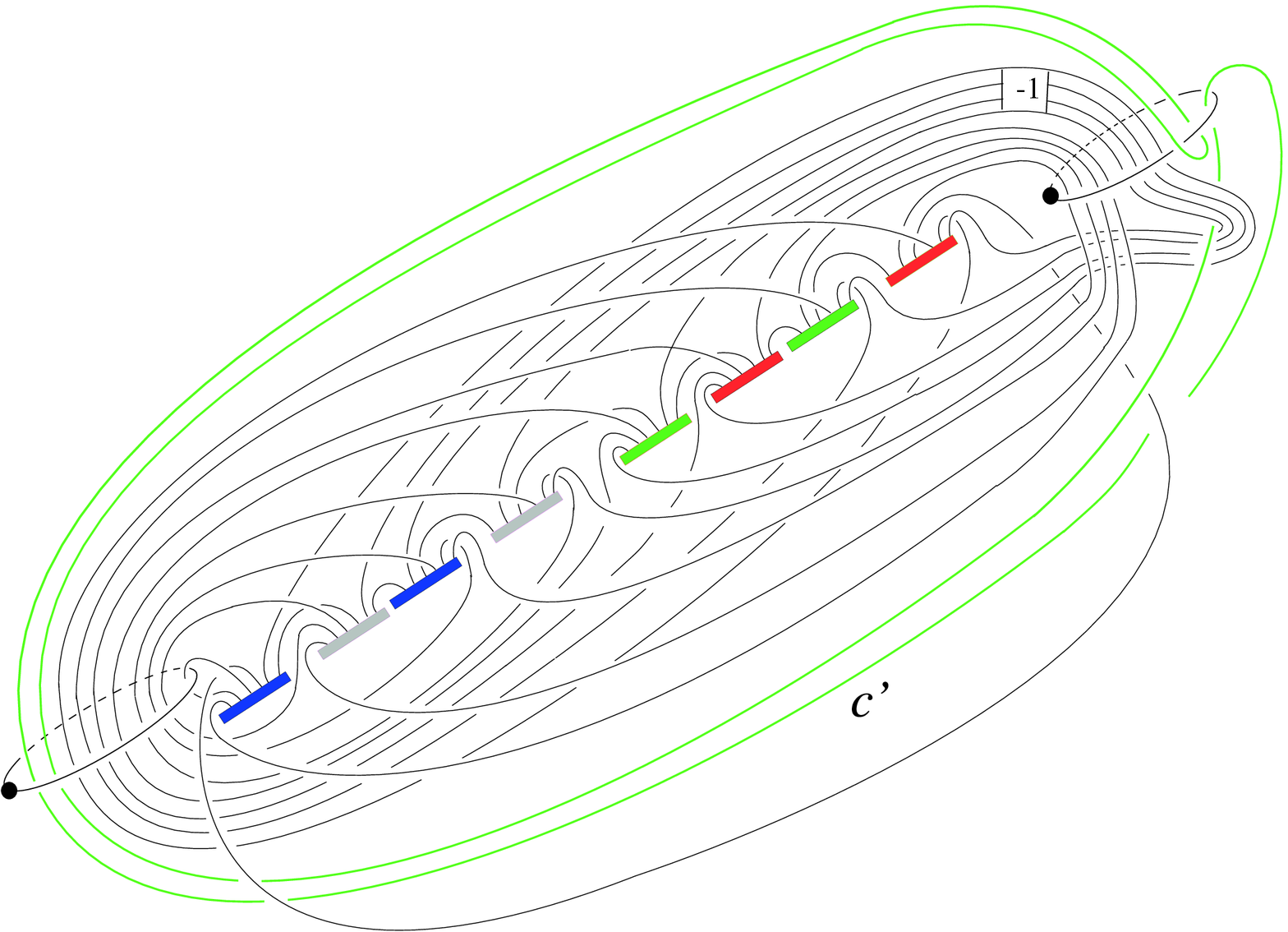}   
\caption{Converting $1$-handle notation} 
\end{center}
\end{figure}

 \begin{figure}[ht]  \begin{center}  
\includegraphics[width=.85\textwidth]{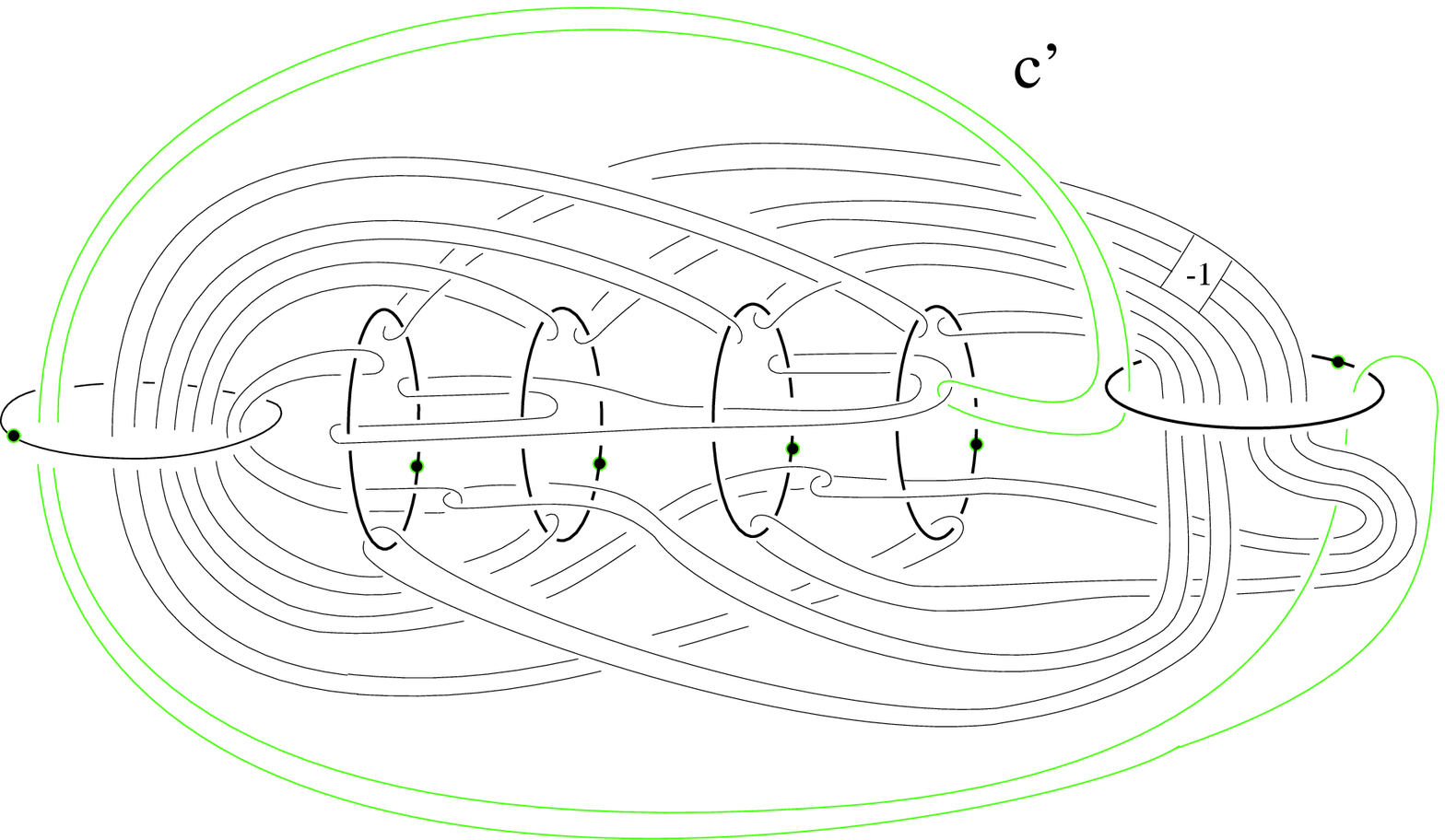}   
\caption{$E_{0}'$} 
\end{center}
\end{figure}

\begin{figure}[ht]  \begin{center}  
\includegraphics[width=.8\textwidth]{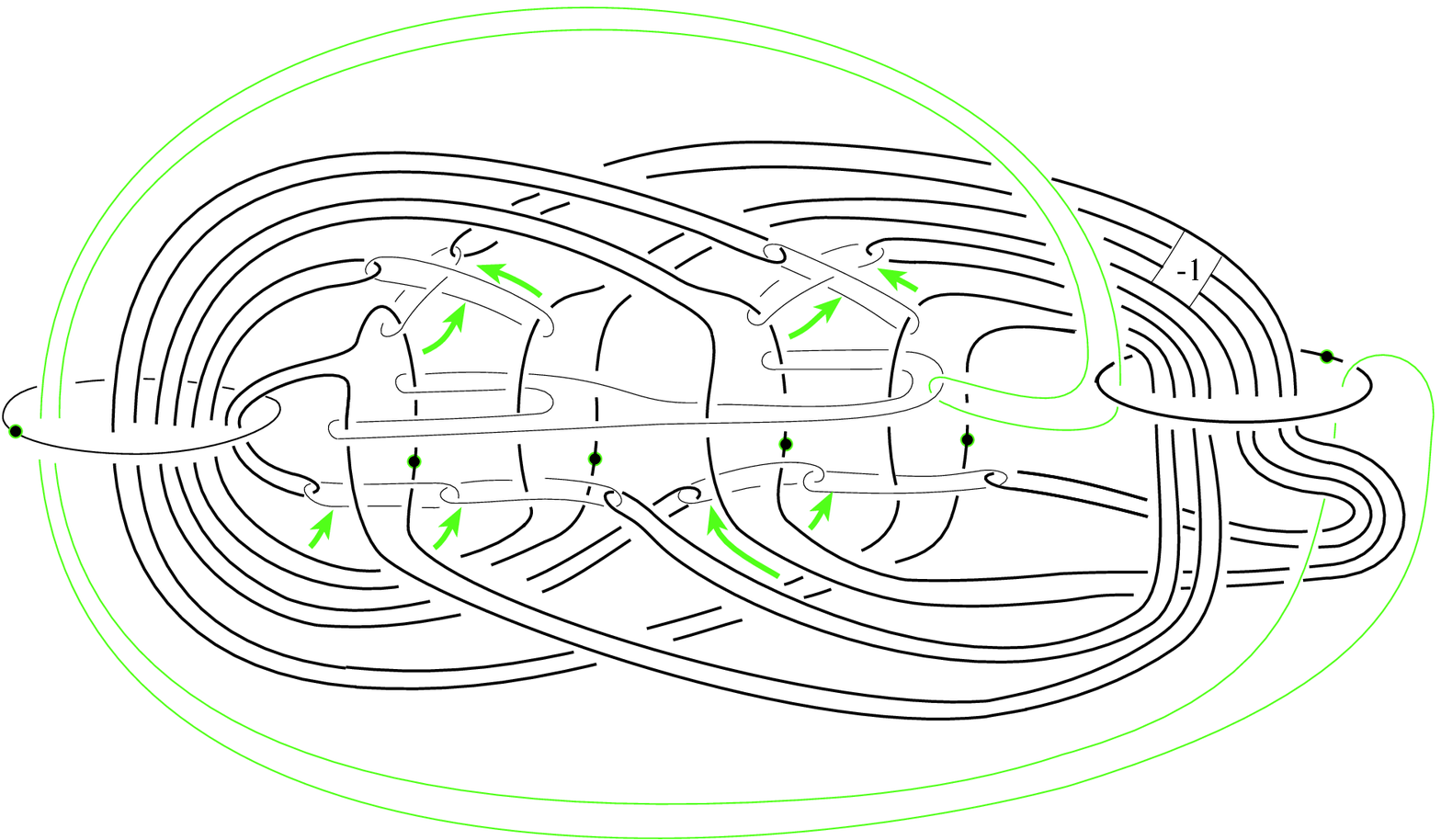}   
\caption{Surgering inside of $E'_{0}$} 
\end{center}
\end{figure}

\newpage

\begin{figure}[ht]  \begin{center}  
\includegraphics[width=.8\textwidth]{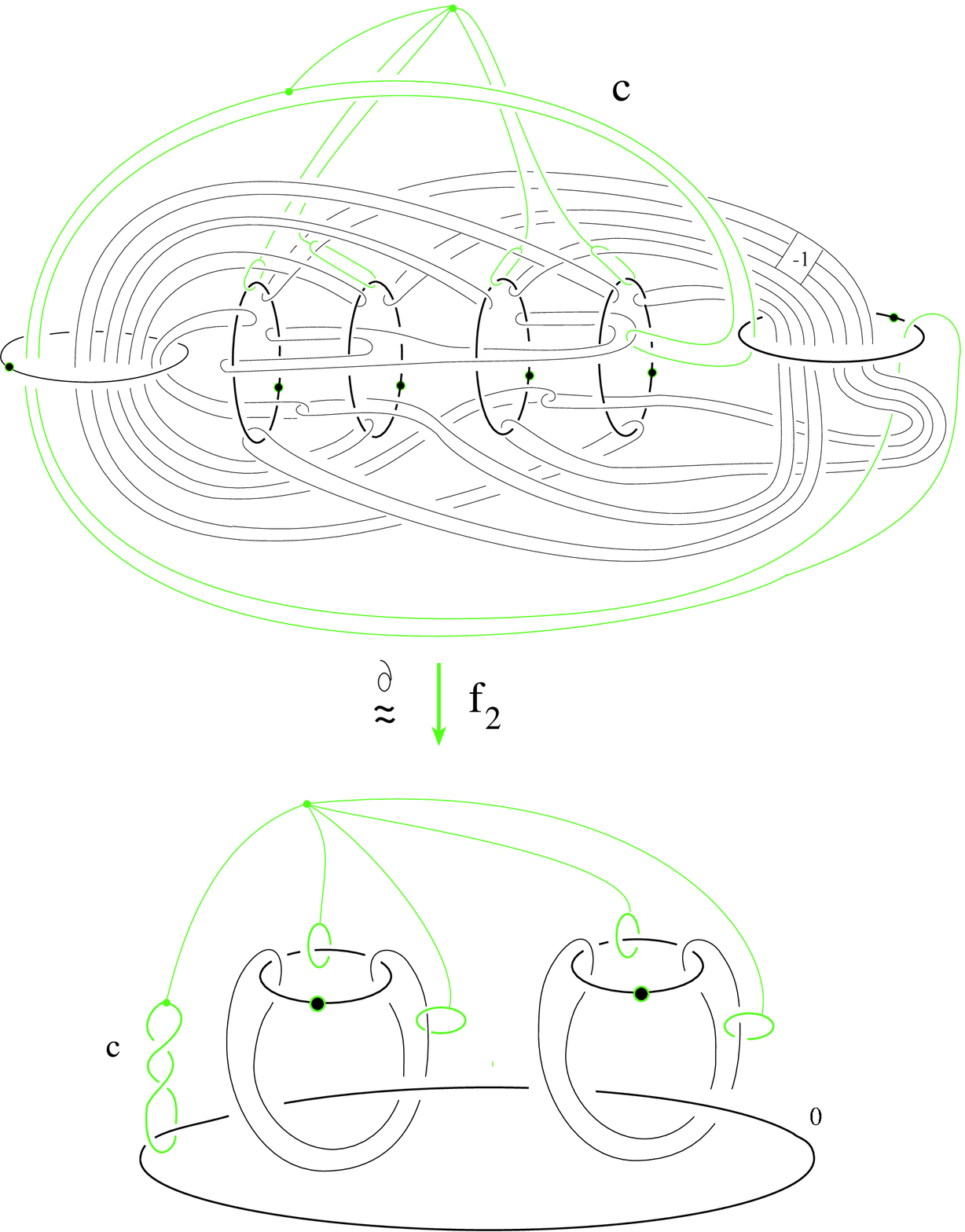}   
\caption{Diffeomorphism $\partial E_{0}'\approx \Sigma_{2}\times S^1$ made concrete} 
\end{center}
\end{figure}

\begin{figure}[ht]  \begin{center}  
\includegraphics[width=.8\textwidth]{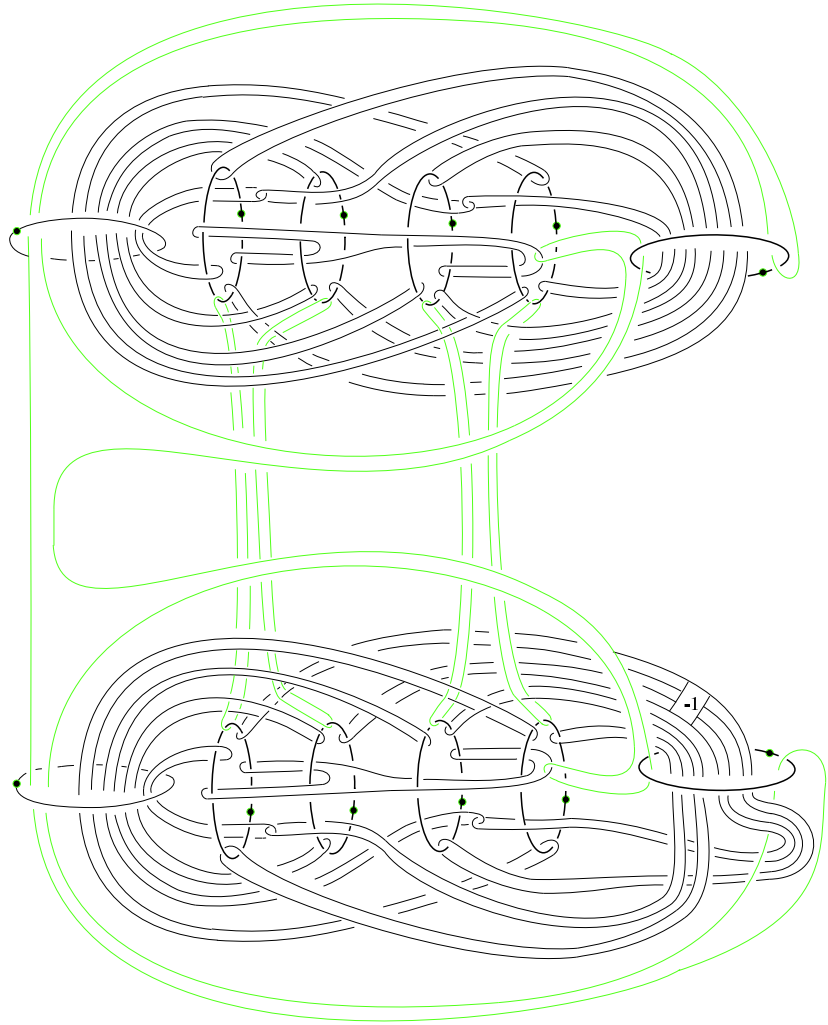}   
\caption{The CaCiMe surface M} 
\end{center}
\end{figure}

\end{document}